\documentclass[11pt]{article}

\usepackage{graphicx,color} 
\usepackage{amsmath,amssymb,latexsym}
\usepackage{mathtools}
\usepackage{xspace}
\usepackage{tikz}
\usepackage[T1]{fontenc}
\usepackage{cite}

\usepackage{hyperref}

\hypersetup{colorlinks=true}

\hypersetup{colorlinks=true, linkcolor=blue, citecolor=blue,urlcolor=blue}

\newtheorem{theorem}{Theorem}[section]
\newtheorem{definition}[theorem]{Definition}
\newtheorem{lemma}[theorem]{Lemma}
\newtheorem{proposition}[theorem]{Proposition}

\newtheorem{conjecture}[theorem]{Conjecture}
\newtheorem{corollary}[theorem]{Corollary}

\newcommand{\proof}{\noindent{\bf Proof.\ }}
\newcommand{\qed}{\hfill $\square$ \bigskip}
\newcommand{\cp}{\,\square\,}

\newcommand{\diam}{{\rm diam}}

\newcommand{\mut}{\mu_{\rm t}}
\newcommand{\mud}{\mu_{\rm d}}

\newcommand{\muo}{\mu_{\rm o}}

\newcommand{\ex}{{\rm ex}}

\title{Mutual-visibility problems on graphs of diameter two}

\author{
Serafino Cicerone $^{a}$\thanks{Email: \texttt{serafino.cicerone@univaq.it}}
\and
Gabriele {Di Stefano} $^{a}$\thanks{Email: \texttt{gabriele.distefano@univaq.it}}
\and
Sandi Klav\v zar $^{b,c,d}$\thanks{Email: \texttt{sandi.klavzar@fmf.uni-lj.si}}
\and Ismael G. Yero $^{e}$\thanks{Email: \texttt{ismael.gonzalez@uca.es}}
}

\begin{document}

\maketitle

\begin{center}
	$^a$ Department of Information Engineering, Computer Science, and Mathematics,
	     University of L'Aquila, Italy \\
	\medskip

	$^b$ Faculty of Mathematics and Physics, University of Ljubljana, Slovenia\\
	\medskip
	
	$^c$ Institute of Mathematics, Physics and Mechanics, Ljubljana, Slovenia\\
	\medskip
	
	$^d$ Faculty of Natural Sciences and Mathematics, University of Maribor, Slovenia\\
	\medskip
	
	$^e$ Departamento de Matem\'aticas, Universidad de C\'adiz, Algeciras Campus, Spain \\

\end{center}

\begin{abstract}
The mutual-visibility problem in a graph $G$ asks for the cardinality of a largest set of vertices $S\subseteq V(G)$ so that for any two vertices $x,y\in S$ there is a shortest $x,y$-path $P$ so that all internal vertices of $P$ are not in $S$. This is also said as $x,y$ are visible with respect to $S$, or $S$-visible for short. Variations of this problem are known, based on the extension of the visibility property of vertices that are in and/or outside $S$. Such variations are called total, outer and dual mutual-visibility problems. This work is focused on studying the corresponding four visibility parameters in graphs of diameter two, throughout showing bounds and/or closed formulae for these parameters. 

The mutual-visibility problem in the Cartesian product of two complete graphs is equivalent to (an instance of) the celebrated Zarankievicz's problem. Here we study the dual and outer mutual-visibility problem for the Cartesian product of two complete graphs and all the mutual-visibility problems for the direct product of such graphs as well. We also study all the mutual-visibility problems for the line graphs of complete and complete bipartite graphs. As a consequence of this study, we present several relationships between the mentioned problems and some instances of the classical Tur\'an problem. Moreover, we study the visibility problems for cographs and several non-trivial diameter-two graphs of minimum size. 
\end{abstract}


\noindent
AMS Subj.\ Class.\ (2020):  05C12, 05C38, 05C69, 05C76

\section{Introduction}

The mutual-visibility problem in graphs has recently appeared in \cite{DiStefano-2022}, and has remarkably attracted the attention of several investigations, which can be seen in the series of articles \cite{Bresar,Bujtas,Cicerone-geo,Cicerone-2023+, Cicerone-2023+a,Cicerone-hered,variety-2023,Cicerone-2023, Cicerone-2022+,kuziak-2023,tian-2023+}. Some reasons of such interest might come from the following facts.
\begin{itemize}
  \item The problem has some origin in a computer science application related to situations arising in a framework of mobile entities in a network. That is, nodes of a network having some ``mutual-visibility'' properties can be seen as entities of a network requiring to communicate between themselves in a somehow confidential or private way. Namely, satisfying that for any exchanged information there should be a channel which does not pass through other entities. For some of these applied researches see for instance \cite{aljohani-2018a,bhagat-2020,Cicerone-2023+,diluna-2017,poudel-2021}.
  \item A close relationship that exists between the mutual-visibility problem and the general position problem~\cite{manuel-2018, ullas-2016}, which is also a distance related topic of high interest in the last recent years~\cite{klavzar-2023, klavzar-2022, klavzar-2021, korze-2023, patkos-2020, tian-2021, yao-2022}, see also~\cite{klavzar-tan-2023, manuel-2022} for the edge version of the general position problem in graphs. 

  \item The connections that have appeared between the mutual-visibility problem with some classical topics in combinatorics. For instance, while studying the mutual-visibility problem in the Cartesian product of complete graphs, it has been noted that solving such a problem turns out to be equivalent to solve an instance of the well-known Zarankievicz's problem (see \cite{Cicerone-2023}). Relatively similar to this, while considering the lower version of this problem, a closed relationship with a classical Bollob\'as-Wessel theorem was proved (see \cite{Bresar}). Also, for the case of the total variant of the mutual-visibility, and the same families of graphs, it has been noted that it can be reformulated as a Tur\'an-type problem on hypergraphs (see \cite{Bujtas}).
  \item The standard mutual-visibility problem can be (and sometimes even needs to be) modified in several directions in order to consider different visibility situations. For instance, while studying the mutual-visibility problem in general Cartesian product graphs (see \cite{Cicerone-2023}), the notion of independent mutual-visibility was naturally required, thus defined, and their first basic properties identified. In the article \cite{Cicerone-2022+}, a total version of the mutual-visibility problem was needed, in order to study the strong product of graphs. This total notion was also a first step into the work \cite{variety-2023}, where this total version was further studied, together with two ``partially'' total ones that were introduced in order to close all the possible ``visibility'' situations that might exist between the elements of a graph.
\end{itemize}

In a formal way, given a connected graph $G$ and a set of vertices $X\subseteq V(G)$, two vertices $x,y\in V(G)$ are called to be $X$-\emph{visible} if there is a shortest $x,y$-path (also called geodesic) whose interior vertices do not belong to $X$. With this idea in mind, for a given set $X\subseteq V(G)$ of a connected graph $G$, the following definitions are known from \cite{variety-2023}.
\begin{itemize}
  \item \emph{Mutual-visibility set}: if any two vertices of $X$ are $X$-visible.
  \item \emph{Outer mutual-visibility set}: if any two vertices $x,y\in X$ and any two vertices $x\in X$ and $y\in \overline{X}$ are $X$-visible.
  \item \emph{Dual mutual-visibility set}: if any two vertices $x,y\in X$ and any two vertices $x,y\in \overline{X}$ are $X$-visible.
  \item \emph{Total mutual-visibility set}: if any two vertices $x,y\in V(G)$ are $X$-visible.
\end{itemize}
Regarding such graph structures, the following parameters are defined as the cardinalities of the largest (respectively) mutual-visibility sets from the above ones.

\medskip
\begin{center}
\begin{tabular}{|c|c|c|c|}
  \hline
  mutual-visibility number & $\mu(G)$ & dual mutual-visibility number & $\mud(G)$ \\ \hline
  outer mutual-visibility number & $\muo(G)$ & total mutual-visibility number & $\mut(G)$ \\
  \hline
\end{tabular}
\end{center}
If $\tau\in \{\mu, \mud, \muo, \mut\}$, then we say that $X\subseteq V(G)$ in a $\tau$-set if $|X| = \tau(G)$. 

In the present investigation, we are focused on giving some contributions on these four mutual-visibility parameters on graphs of diameter two. Some motivations for this specific study are coming from already established results on such class of graphs. For instance, as already mentioned, the mutual-visibility problem in the Cartesian product of complete graphs is proved to be equivalent to solve an instance of the well-known Zarankievicz's problem (see \cite{Cicerone-2023}), and such Cartesian products are of diameter two. In this sense, we continue this research direction on graphs of diameter two, which  will indeed show that this problem, and the related variations, remain challenging while considering diameter-two graphs in general.

In the remaining of this section we give some preliminary terminologies and notations that shall be used throughout our exposition. In Section \ref{sec:hamming} we consider the Cartesian and direct products of complete graphs. Specifically, we give formulas for the dual and outer mutual-visibility numbers of the Cartesian product, which fulfills the existing gap for the visibility numbers of such graphs. These results allow us, among other things, to answer negatively a question in the literature regarding the relationship between the mutual-visibility number and the outer mutual-visibility number. We also compute all the mutual-visibility numbers of the direct products, showing that all of them achieve the same value.
Section \ref{sec:line-graphs} focuses on the line graphs of complete and complete bipartite graphs. Through this study, we give several relationships between the mutual-visibility problems and some instances of the classical Tur\'an problem. Among them, we for instance show that the mutual-visibility number of the line graph of complete graphs equals the number of edges of the Tur\'an graph $T(n,3)$, and that the total mutual-visibility number of such graphs equals the number of edges of the Tur\'an related graph $\ex(n;C_4)$. Connections between the mutual-visibility problem on the line graphs of complete bipartite graphs and the Zarankiewicz problem are also given in Section \ref{sec:line-graphs}. Next, in Section \ref{sec:cographs} we consider the class of cographs, by studying those graphs $G$ that have values in their mutual-visibility numbers equal to at least the order of $G$ minus one. Section \ref{sec:of-minimum-size} is focused on non-trivial diameter-two graphs of minimum size. That is, we compute the values of the mutual-visibility parameters of the graphs belonging to this class. Finally, Section \ref{sec:conclusion} gives some concluding remarks together with some future research lines that can be of interest as a continuation of this work. 

\subsection{Preliminaries}

All graphs considered in this paper are finite and simple. The \emph{distance} $d_G(u,v)$ between vertices $u$ and $v$ of a graph $G$ is the length of a shortest $u,v$-path. 
The degree of a vertex $v$ in $G$ is denoted as $\deg_G(v)$. The {\em girth}, $g(G)$, of a graph $G$ is the length of a shortest cycle of $G$. If $G$ is a forest, then we set $g(G) = \infty$. For an integer $k\ge 1$, we shall write $[k]=\{1,\dots,k\}$. The order and the size of $G$ will be respectively denoted by $n(G)$ and $m(G)$. A vertex of $G$ is {\em universal} if it is adjacent to all the other vertices of $G$. 

A {\em cograph} is a graph which contains no induced path on four vertices. Cographs can be characterized in many different ways, see~\cite{corneil-1981}. For instance, cographs are precisely the graphs that can be obtained from $K_1$ by means of the disjoint union and join of graphs.

The \emph{Cartesian product} $G\cp H$ and the \emph{direct product} $G\times H$ of graphs $G$ and $H$ both have the vertex set $V(G)\times V(H)$. In $G\cp H$, vertices $(g,h)$ and $(g',h')$ are adjacent if either $g=g'$ and $hh'\in E(H)$, or $h=h'$ and $gg'\in E(G)$. In $G\times H$, vertices $(g,h)$ and $(g',h')$ are adjacent if $gg'\in E(G)$ and $hh'\in E(H)$. In each of the two products, if $h\in V(H)$, then the set of vertices $\{(g,h):\ g\in V(G)\}$ forms a {\em $G$-layer} which is denoted by $G^h$. For a given $g\in V(G)$, the {\em $H$-layer} $^gH$ is defined analogously. Note that in $G\cp H$, layers induce subgraphs isomorphic to $G$ resp.\ $H$, while in $G\times H$ layers induce edgeless graphs. 

The \emph{union} $G\cup H$ of $G$ and $H$ is the graph $G'$  with vertex set $V(G') = V(G)\cup V(H)$ and edge set $E(G') = E(G)\cup E(H)$. 
The \emph{join} $G + H$ of $G$ and $H$ is the graph $G'$ with vertex set $V(G') = V(G)\cup V(H)$ and edge set $E(G') = E(G)\cup E(H) \cup \{uv~:~ u\in V(G), u\in V(H)\}$. 

\smallskip
We conclude these preliminaries with the following result. 

\begin{lemma}
\label{lem:girth5}
If $G$ is a connected graph of order at least $3$ and with $g(G)\ge 5$, then an outer mutual-visibility set is an independent set. 
\end{lemma}

\proof
Let $X$ be an outer mutual-visibility set of $G$ and suppose that $X$ contains vertices $x$ and $y$ such that $xy\in E(G)$. Since $G$ has at least three vertices and is connected, we may assume that $z$ is a neighbor of $y$ different from $x$. As $g(G)\ge 5$, we must have $z\in X$, for otherwise $x\in X$ and $z\notin X$ are not $X$-visible. But then $x, z\in X$ are not $X$-visible, a contradiction. 
\qed

\section{Products of two complete graphs}
\label{sec:hamming}

In this section we consider the four invariants of interest on Cartesian and direct products of two complete graphs. The invariants $\mu$ and $\mu_t$ were already considered, here we add formulas for $\mud$ and $\muo$. Using these results we are able to answer in negative the question from~\cite{variety-2023} whether $\mu(G) \le 2\muo(G)$ holds for any graph $G$. On the other hand, the direct product of (complete) graphs was not yet considered in this context, in this section we prove the formula for it which is the same for all the four invariants. 

In~\cite{Cicerone-2023} it is proved that if $m, n\ge 2$, then $\mu(K_m\cp K_n) = z(m,n;2,2)$, where $z(m, n; 2, 2)$ is the maximum number of $1$s that an $m\times n$ binary matrix can have, provided that it contains no $2\times 2$ submatrix of $1$s. To determine $z(m, n; 2, 2)$ is a notorious open (instance of)  Zarankievicz's problem. When $n$ is sufficiently large, the value $z(n, n; 2, 2)$ can be bounded as follows~\cite{brown-1966, erdos-1966}:
$$n^{3/2}-n^{4/3} \le z(n,n;2,2) \le \frac{1}{2} n(1+ \sqrt{4n-3})\,,$$ 
which demonstrates that the growth is faster than linear. 

For the total mutual-visibility it was proved in~\cite{tian-2023+} that if $n, m\ge 2$, then 
$$\mut(K_n\cp K_m) = \max \{n,m\}\,.$$
For the dual and the outer mutual-visibility number, we have the following related respective results. 

\begin{theorem}
\label{thm:hamming-outer-dual}
If $n, m\ge 3$, then 
\begin{enumerate}
    \item[{\rm(i)}] $\mud(K_n \cp K_m) = n + m - 1$,
    \item[{\rm(ii)}] $\muo(K_n \cp K_m) = n + m - 2$.
\end{enumerate}
\end{theorem}

\proof
Let $V(K_k) = [k]$, so that $V(K_n \cp K_m) = \{(i,j):\ i\in [n], j\in [m]\}$. For the rest of the proof set $G = K_n \cp K_m$.

(i) It is straightforward to check that the set $\{(i,1):\ i\in [n]\}\cup \{(1,j):\ j\in [m]\}$ is a dual mutual-visibility set of cardinality $n+m-1$, thus $\mud(G)\ge n + m - 1$. 

\begin{figure}[ht!]
\begin{center}
\begin{tikzpicture}[scale=0.4,style=thick,x=1cm,y=1cm]
\def\vr{6pt}

\begin{scope} 
\coordinate(x) at (2,2);
\coordinate(y) at (4,2);
\coordinate(x') at (4,4);
\coordinate(a) at (2,4);
\draw (0,-0.4) -- (0,10.4);
\draw (2,-0.4) -- (2,10.4);
\draw (4,-0.4) -- (4,10.4);
\draw (6,-0.4) -- (6,10.4);
\draw (8,-0.4) -- (8,10.4);
\draw (10,-0.4) -- (10,10.4);
\draw (-0.4,0) -- (10.4,0);
\draw (-0.4,2) -- (10.4,2);
\draw (-0.4,4) -- (10.4,4);
\draw (-0.4,6) -- (10.4,6);
\draw (-0.4,8) -- (10.4,8);
\draw (-0.4,10) -- (10.4,10);
\draw(x)[fill=black] circle(\vr);
\draw(y)[fill=black] circle(\vr);
\draw(x')[fill=black] circle(\vr);
\draw(a)[fill=white] circle(\vr);
\node at (5,-3) {(a)};
\node at (2,-1) {$i$};
\node at (4,-1) {$i+1$};
\node at (-1,2) {$j$};
\node at (-1.5,4) {$j+1$};
\node at (2.5,1.5) {$x$};
\node at (4.5,1.5) {$y$};
\node at (4.5,3.5) {$x'$};
\end{scope}
		
\begin{scope}[xshift=15cm, yshift=0cm] 
\coordinate(x) at (2,2);
\coordinate(y) at (4,2);
\coordinate(x') at (4,4);
\coordinate(a) at (2,4);
\coordinate(z) at (2,8);
\coordinate(b) at (4,8);
\draw (0,-0.4) -- (0,10.4);
\draw (2,-0.4) -- (2,10.4);
\draw (4,-0.4) -- (4,10.4);
\draw (6,-0.4) -- (6,10.4);
\draw (8,-0.4) -- (8,10.4);
\draw (10,-0.4) -- (10,10.4);
\draw (-0.4,0) -- (10.4,0);
\draw (-0.4,2) -- (10.4,2);
\draw (-0.4,4) -- (10.4,4);
\draw (-0.4,6) -- (10.4,6);
\draw (-0.4,8) -- (10.4,8);
\draw (-0.4,10) -- (10.4,10);
\draw(x)[fill=black] circle(\vr);
\draw(y)[fill=black] circle(\vr);
\draw(x')[fill=black] circle(\vr);
\draw(a)[fill=white] circle(\vr);
\draw(z)[fill=black] circle(\vr);
\draw(b)[fill=white] circle(\vr);
\node at (5,-3) {(b)};
\node at (2,-1) {$i$};
\node at (4,-1) {$i+1$};
\node at (-1,2) {$j$};
\node at (-1.5,4) {$j+1$};
\node at (-1.1,8) {$j''$};
\node at (2.5,1.5) {$x$};
\node at (4.5,1.5) {$y$};
\node at (2.5,7.5) {$z$};
\node at (4.5,3.5) {$x'$};
\end{scope}

\begin{scope}[xshift=0cm, yshift=-15cm] 
\coordinate(x) at (2,2);
\coordinate(y) at (4,2);
\coordinate(x') at (4,4);
\coordinate(a) at (2,4);
\coordinate(z) at (8,8);
\coordinate(c) at (8,4);
\coordinate(d) at (8,2);
\draw (0,-0.4) -- (0,10.4);
\draw (2,-0.4) -- (2,10.4);
\draw (4,-0.4) -- (4,10.4);
\draw (6,-0.4) -- (6,10.4);
\draw (8,-0.4) -- (8,10.4);
\draw (10,-0.4) -- (10,10.4);
\draw (-0.4,0) -- (10.4,0);
\draw (-0.4,2) -- (10.4,2);
\draw (-0.4,4) -- (10.4,4);
\draw (-0.4,6) -- (10.4,6);
\draw (-0.4,8) -- (10.4,8);
\draw (-0.4,10) -- (10.4,10);
\draw(x)[fill=black] circle(\vr);
\draw(y)[fill=black] circle(\vr);
\draw(x')[fill=black] circle(\vr);
\draw(z)[fill=black] circle(\vr);
\draw(c)[fill=black] circle(\vr);
\draw(d)[fill=black] circle(\vr);
\draw(a)[fill=white] circle(\vr);
\node at (5,-3) {(c)};
\node at (2,-1) {$i$};
\node at (4,-1) {$i+1$};
\node at (8,-1) {$i''$};
\node at (-1,2) {$j$};
\node at (-1.5,4) {$j+1$};
\node at (-1.1,8) {$j''$};
\node at (2.5,1.5) {$x$};
\node at (4.5,1.5) {$y$};
\node at (4.5,3.5) {$x'$};
\node at (7.5,7.5) {$z$};
\end{scope}

\begin{scope}[xshift=15cm, yshift=-15cm] 
\coordinate(x) at (2,2);
\coordinate(y) at (4,2);
\coordinate(x') at (4,4);
\coordinate(a) at (2,4);
\coordinate(z) at (8,8);
\coordinate(c) at (8,4);
\coordinate(d) at (4,8);
\coordinate(e) at (2,8);
\coordinate(f) at (8,2);
\draw (0,-0.4) -- (0,10.4);
\draw (2,-0.4) -- (2,10.4);
\draw (4,-0.4) -- (4,10.4);
\draw (6,-0.4) -- (6,10.4);
\draw (8,-0.4) -- (8,10.4);
\draw (10,-0.4) -- (10,10.4);
\draw (-0.4,0) -- (10.4,0);
\draw (-0.4,2) -- (10.4,2);
\draw (-0.4,4) -- (10.4,4);
\draw (-0.4,6) -- (10.4,6);
\draw (-0.4,8) -- (10.4,8);
\draw (-0.4,10) -- (10.4,10);
\draw(x)[fill=black] circle(\vr);
\draw(y)[fill=black] circle(\vr);
\draw(x')[fill=black] circle(\vr);
\draw(z)[fill=black] circle(\vr);
\draw(c)[fill=white] circle(\vr);
\draw(d)[fill=black] circle(\vr);
\draw(a)[fill=white] circle(\vr);
\draw(e)[fill=white] circle(\vr);
\draw(f)[fill=white] circle(\vr);
\node at (5,-3) {(d)};
\node at (2,-1) {$i$};
\node at (4,-1) {$i+1$};
\node at (8,-1) {$i''$};
\node at (-1,2) {$j$};
\node at (-1.5,4) {$j+1$};
\node at (-1.1,8) {$j''$};
\node at (2.5,1.5) {$x$};
\node at (4.5,1.5) {$y$};
\node at (4.5,3.5) {$x'$};
\node at (7.5,7.5) {$z$};
\end{scope}
\end{tikzpicture}
\caption{Cases from the proof of Theorem~\ref{thm:hamming-outer-dual}(i)}
	\label{fig:cases}
\end{center}
\end{figure}
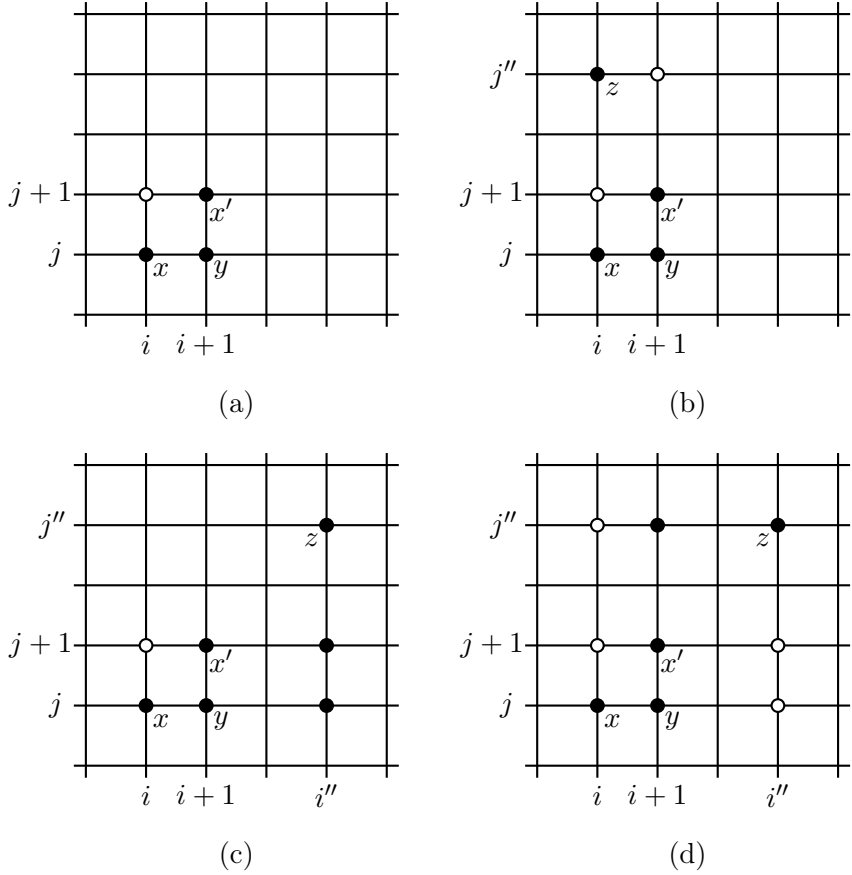

To prove the reverse inequality, suppose for a purpose of contradiction that there exists a dual mutual-visibility set $X$ of $G$ of cardinality $n+m$. Let $x = (i,j)$ be an arbitrary vertex from $X$. Since the union of the two layers containing $x$ contains $n+m-1$ vertices, there exists a vertex $x' = (i', j')\in X$, where $i\ne i'$ and $j\ne j'$. By the symmetry of $G$ we may without loss of generality assume that $x'=(i+1,j+1)$. (Here and later below, indices are computed modulo $n$ and $m$ if necessary.) As $X$ is a dual mutual-visibility set, at least one of the vertices $(i+1,j)$ and $(i,j+1)$ must belong to $X$, for otherwise they are not $X$-visible. We may without loss of generality assume that $y = (i+1,j)\in X$. Then $(i,j+1)\notin X$. In Fig.~\ref{fig:cases}(a) the situation so far is schematically presented, where we use the convention that the vertices from $X$ are shown in black, and the vertices not from $X$ in white. 

Consider now the vertex $y=(i+1,j)$. By the above argument used for $x$, there exists a vertex $z\in X$ which does not lie in the union of the two layers containing $y$. Assume first that $z=(i,j'')$, where $j''\ne j$. Then clearly we also have $j''\ne j+1$. Now we must have $(i+1,j'')\not\in X$, for otherwise $y$ and $z$ are not $X$-visible. But then none of the vertices $(i+1,j'')$ and $(i,j+1)$ lies in $X$ and are not $X$-visible, a contradiction. This situation is shown in Fig.~\ref{fig:cases}(b).

If $z=(i'',j+1)$, then by symmetry, we also arrive to a contradiction. 

Assume now that $z=(i'',j'')$, where $i''\ne i, i+1$ and $j''\ne j, j+1$. Considering the vertices $z=(i'',j'')\in X$ and $x'=(i+1,j+1)\in X$, we see that one of $(i'', j+1)$ and $(i+1, j'')$ belongs to $X$; for otherwise we get a contradiction, since these two vertices are not $X$-visible. Suppose that $(i'', j+1)\in X$, see Fig.~\ref{fig:cases}(c). Consider the vertices $x=(i,j)\in X$, $(i'', j+1)\in X$, and $(i,j+1)\notin X$ to realize that $(i'',j)\in X$. But then $(i'',j)$ and $x'$ are not $X$-visible. This implies that $(i'', j+1)\notin X$, and consequently $(i+1, j'')\in X$, see Fig.~\ref{fig:cases}(d). Now $(i,j'')\notin X$, for otherwise these vertex cannot see $y$. Similarly, $(i'', j)\notin X$, see Fig.~\ref{fig:cases}(d) again. But now the vertices $(i,j'')\notin X$ and $(i'', j)\notin X$, are not $X$-visible. This final contradiction implies that $\mud(G)\le n + m - 1$ which proves (i). 

\medskip
(ii) It is straightforward to verify that the set $\{(i,1):\ 2\le i\le n\}\cup \{(1,j):\ 2\le j\le m\}$ is an outer mutual-visibility set of cardinality $n+m-2$, thus $\muo(G)\ge n + m - 2$. 

Let $X$ be an arbitrary outer mutual-visibility set of $G$. We may assume without loss of generality that $(1,1)\notin X$. Let $x=(i,j)$ be an arbitrary vertex of $X$. Then in one of the layers $(K_n)^j$ and $^i(K_m)$, the vertex $x$ is the only vertex from $X$. Indeed, if we would have $(i,j')\in X$, $j'\ne j$, and $(i',j)\in X$, $i'\ne i$, then no matter whether $(i',j')$ lies in $X$ or not, the vertex $x$ would not see $(i',j')$. We are now going to assign to each vertex $x\in X$ a unique variable as follows. If $x=(i,1)\in X$, then assign to $x$ the variable $a_i$ and if  $x=(1,j)\in X$, then assign to $x$ the variable $b_j$. In addition, if $i, j\ge 2$ and $x=(i,j)\in X$, then in the case that $X \cap V((K_n)^j) = \{x\}$, we assign to $x$ the variable $b_j$, and if $X \cap V(^i(K_m)) = \{x\}$, then we assign to $x$ a variable $a_i$. Note that if, say, $(i,1)\in X$ and $(i,j)\in X$, $j\ne 1$, then $(i,1)$ is assigned $a_i$ and $(i,j)$ is assigned $b_j$. 
Since to each vertex of $X$ we assign a different variable and, having in mind that $(1,1)\notin X$, the variables used are $a_2, \ldots, a_n$ and $b_2, \ldots, b_m$, we have $|X|\le n+m-2$. We conclude that $\muo(G)\le n + m - 2$. 
\qed

By Theorem~\ref{thm:hamming-outer-dual} and the discussion before it, as soon as $n$ and $m$ are not small, we have 
\begin{equation}
\label{eq:cart}    
\mut(K_n \cp K_m) < \muo(K_n \cp K_m) < \mud(K_n \cp K_m) < \mu(K_n \cp K_m)\,.
\end{equation}

In~\cite{variety-2023} a question was posed whether $\mu(G) \le 2\muo(G)$ is true in general. We can now answer this question in negative because $\muo(K_n \cp K_n) = 2n - 2$ and $\mu(K_n\cp K_n) \ge n^{3/2}-n^{4/3}$.

We now turn our attention to the direct product of complete graphs and prove the following result. 

\begin{theorem}
\label{thm:direct-product}
If $n, m\ge 5$, then $\mut(K_n \times K_m) = \mu(K_n \times K_m) = nm - 4$.
\end{theorem}

\proof
Let $V(K_k) = [k]$, so that $V(K_n \times K_m) = \{(i,j):\ i\in [n], j\in [m]\}$. Set $G = K_n \times K_m$ for the rest of the proof. 

We claim first that $\mu(G) \le nm-4$. Suppose on the contrary that there exists a mutual-visibility set $X$ with $|X| = nm-3$. Let $V(G)\setminus X = \{x,y,z\}$. By the symmetry of $G$, it suffices to consider the following four cases. 

Suppose that $x= (i,j)$, $y= (i',j')$, $z= (i'',j'')$, where $|\{i,i',i''\}| = 3$ and  $|\{j,j',j''\}| = 3$. In this case we see that the vertices $(i',j)$ and $(i',j'')$ belong to $X$ but are not $X$-visible. Hence this case is not possible. 

In the second case, suppose that $x= (i,j)$, $y= (i',j)$, $z= (i',j')$, where $i\ne i'$ and $j\ne j'$. Now we infer that the vertices $(i-1,j)$ and $(i-1,j')$ belong to $X$ but are not $X$-visible. 

In the third case, suppose that $x= (i,j)$, $y= (i,j')$, $z= (i',j'')$, where $i\ne i'$ and $|\{j,j',j''\}| = 3$. Now consider the vertices $(i',j)$ and $(i',j')$ which both belong to $X$ but are not $X$-visible. 

In the last case, suppose that $x= (i,j)$, $y= (i',j)$, $z= (i'',j)$, where $|\{i,i',i''\}| = 3$. Now consider two vertices $(k,j)$ and $(k',j)$, where $k$ and $k'$ are selected in such a way that $|\{i,i',i'',k,k'\}| = 5$. (Such values $k$ and $k'$ exist since we have assumed that $n,m \ge 5$.) But now the vertices $(k,j)$ and $(k',j)$ belong to $X$, and they are not $X$-visible.  

We can conclude that no matter how the set $X$ lies in $G$, it cannot form a mutual-visibility set. This proves that $\mu(G) \le nm - 4$.

Let $Y = \{(1,1), (2,2), (3,3), (4,4)\}$ and let $X = V(G) \setminus Y$. We claim that $X$ is a total mutual-visibility set of $G$. Let $x=(i,j)$ and $y=(i',j')$ be arbitrary vertices of $G$. Assume first that $x,y\in X$. If $i\ne i'$ and $j\ne j'$, then $xy\in E(G)$ and there is nothing to prove. Otherwise, $i=i'$ or $j=j'$. Assume without loss of generality that $i=i'$ and let $k\in [4]$ be such that $k\ne i, j, j'$. Then $(i,j)(k,k)\in E(G)$ and $(k,k)(i,j')\in E(G)$, hence $x$ and $y$ are $X$-visible. We proceed similarly in the case when $x\in X$ and $y\in Y$. Finally, if $x,y\in Y$, then $xy\in E(G)$. We have thus demonstrated that $X$ is a total mutual-visibility set of $G$.  

By the above, $\mut(G) \ge mn- 4$. Combining this inequality with the earlier proved inequality $\mu(G) \le nm - 4$, we have 
$$nm - 4  \le \mut(G) \le \mu(G) \le nm - 4\,,$$
hence the equality holds everywhere and we are done. 
\qed

\begin{corollary}
\label{cor:direct-product}
If $n,m\ge 5$, then 
$$\mut(K_n \times K_m) = \muo(K_n \times K_m) = \mud(K_n \times K_m) = \mu(K_n \times K_m)\,.$$    
\end{corollary}

\proof
Combine Theorem~\ref{thm:direct-product} with the facts following directlty from definitions that for any graph $G$ we have $\mut(G)\le \muo(G)\le \mu(G)$ and $\mut(G)\le \mud(G)\le \mu(G)$, cf.~\cite{variety-2023}. 
\qed

Note that Corollary~\ref{cor:direct-product} is in sharp contrast to~\eqref{eq:cart}. 

\section{Line graphs}
\label{sec:line-graphs}

Given a graph $G$, the \emph{line graph} $L(G)$ of $G$ has vertex set $V(L(G))=\{e_{uv} \,:\,uv\in E(G)\}$, and two vertices $e_{uv},e_{u'v'}$ are adjacent in $L(G)$ if and only if the edges $uv,u'v'$ are incident in $G$. From now on, given a set of edges $F\subseteq E(G)$, we set $S_F = \{e_{uv}\in V(L(G)): \ uv\in F\}$. Also, by $G_{F}$ we represent the subgraph of $G$ whose edges are those ones in $F$ and vertices are those from the edges of $F$.

In this section we focus on the line graphs of complete graphs and of complete bipartite graphs. Notice that if $n\ge 4$, then $\diam(L(K_n)) = 2$, and if $m, n\ge 2$, then $\diam(L(K_{m,n})) = 2$. More generally, if $\diam(G) \le 2$, then $\diam(L(G))\le 3$. We begin with a characterization of mutual visibility sets in line graphs $L(G)$ for graphs $G$ with $\diam(G) = 2$. 

\begin{lemma}
\label{lem:MVS-line}
Let $G$ be a graph of diameter $2$ and $F\subseteq E(G)$. Then $S_F\subseteq V(L(G))$ is a mutual-visibility set of $L(G)$ if and only if for any two independent edges $uv,u'v'\in F$ one of the following conditions is satisfied.
\begin{itemize}
  \item[{\rm(i)}] There is an edge $xy\notin F$ incident with both $uv$ and $u'v'$, or
  \item[{\rm(ii)}] $d_{L(G)}(e_{uv},e_{u'v'}) = 3$ and there is a vertex $z\in V(G)$ adjacent to (w.l.g.) $u$ and $u'$ in $G$, such that $uz, u'z\notin F$.
\end{itemize}
\end{lemma}

\proof
$(\Rightarrow)$ 
Assume $S_F$ is a mutual-visibility set of $L(G)$, and let $uv,u'v'\in F$ be two independent edges. Since $\diam(G) = 2$, we have $d_{L(G)}(e_{uv},e_{u'v'}) \in \{2,3\}$. If $d_{L(G)}(e_{uv},e_{u'v'}) = 2$, then because $S_F$ is a mutual-visibility set, there exists a vertex $e_{xy}\notin S_F$ such that $e_{uv}e_{xy}e_{u'v'}$ is a geodesic. Thus (i) holds as $xy$ is incident to $uv$ and to $u'v'$. If $d_{L(G)}(e_{uv},e_{u'v'}) = 3$, then there must be a vertex $z\in V(G)$ adjacent to (w.l.g.) $u$ and $u'$ in $G$. Clearly, since $S_F$ is a mutual-visibility set, there must be such vertex $z$ with $e_{uz},e_{u'z}\notin S_F$, hence (ii) holds. 

$(\Leftarrow)$ We need to show that any two vertices $e_{uv},e_{u'v'}\in S_F$ are $S_F$-visible in $L(G)$. There is nothing to prove if $e_{uv}e_{u'v'}\in E(L(G))$. Assume that $d_{L(G)}(e_{uv},e_{u'v'}) = 2$. Then (ii) does not apply, hence by (i) there is an edge $xy$ incident with both $uv,u'v'$ such that $e_{xy}\notin S_F$. Thus, $e_{uv}e_{xy}e_{u'v'}$ is a geodesic whose internal vertices are not in $S_F$, and so, $e_{uv},e_{u'v'}\in S_F$ are $S_F$-visible. Assume next that $d_{L(G)}(e_{uv},e_{u'v'}) = 3$. Then (ii) applies, so that there is a vertex $z\in V(G)$ adjacent to $u$ and $u'$ in $G$, such that $e_{uz},e_{u'z}\notin S_F$. Hence, $e_{uv}e_{uz},e_{u'z}e_{u'v'}$ is a geodesic in $L(G)$ whose interior vertices are not in $S_F$. This $e_{uv},e_{u'v'}$ are $S_F$-visible in $L(G)$ in this case as well. Since $\diam(L(G))\le 3$ we are done. 
\qed

Lemma~\ref{lem:MVS-line} reduces the verification whether a set of vertices of $L(G)$, where $\diam(G)\le 2$, is a  mutual-visibility set to the search for the set of edges of largest cardinality in $G$ satisfying the conditions of the lemma. This can be interpreted as an instance of a Tur\'an-type problem. The first striking example of this claim is the following result. For its statement recall that the {\em Tur\'an graph} $T(n,r)$ is a complete $r$-partite graph of order $n$ in which sizes of the $r$ parts are as equal as possible. 

\begin{theorem}
\label{th:mu-Turan}
Let $n\ge 3$ be an integer and $F\subseteq E(K_n)$. Then $S_F\subseteq V(L(K_n))$ is a $\mu$-set of $L(K_n)$ if and only if $(K_n)_F \cong T(n,3)$.
\end{theorem}

\proof
For $n=3$ we have $L(K_3) = K_3$ and the assertion is clear. Suppose in the rest that $n\ge 4$. Then $\diam(L(K_n)) = 2$ and thus Lemma~\ref{lem:MVS-line} implies that if $S_F$ is a mutual-visibility set of $L(K_n)$, then for any two independent edges $uv,u'v'\in F$ there is an edge $xy$ incident with both $uv,u'v'$ such that $xy\notin F$. This can be equivalently reformulated by saying that $S_F$ is a mutual-visibility set of $L(K_n)$ if and only if $(K_n)_F$ does not contain a $K_4$. Tur\'an's theorem (see~\cite[Theorem 11.1.3]{west-2021}) completes the argument. 
\qed

Since the Tur\'an graph $T(n,r)$ has $(1-\frac{1}{r}+o(1))\,\frac{n^2}{2}$ edges, we deduce the following consequence of Theorem~\ref{th:mu-Turan}.

\begin{corollary}
If $n\ge 3$, then $\mu(L(K_n))=(\frac{2}{3}+o(1))\,\frac{n^2}{2}$.
\end{corollary}

Now, with respect to the remaining mutual-visibility parameters of the graph $L(K_n)$, we note the following facts. If $F$ is a set of edges of $K_n$, then the corresponding set $S_F$ in $L(K_n)$ has the (total, outer or dual) mutual-visibility properties based on the existence of certain structures obtained from pairs of not incident edges from $E(K_n)$, $F$, or $E(K_n)\setminus F$. Recall that a pair of not incident edges $uv,u'v'\in E(K_n)$ are $S_F$-visible in $L(K_n)$ whenever there is an edge $xy\notin F$ such that (w.l.g.) $x=u$ and $y=u'$. These facts, the definitions of (total, outer or dual) mutual-visibility sets and the structure of $L(K_n)$ allow to readily observe the following result, whose proof is rather simple and left to the reader.

\begin{lemma}
\label{lem:TOD-L-K_n}
Let $n\ge 3$ be an integer and let $F\subseteq E(K_n)$. Then,
\begin{itemize}
  \item[{\rm (i)}] $S_F$ is a total mutual-visibility set of $L(K_n)$ if and only if for any two not incident edges $uv,u'v'\in E(K_n)$ the subgraph induced by $u,v,u',v'$ has at least one edge not in $F$ different from $uv$ and $u'v'$.
  \item[{\rm (ii)}] $S_F$ is a dual mutual-visibility set of $L(K_n)$ if and only if 
  \begin{itemize}
    \item for any two not incident edges $uv,u'v'\in E(K_n)\setminus F$ the subgraph induced by $u,v,u',v'$ has at least one edge not in $F$ different from $uv$ and $u'v'$, and
    \item for any two not incident edges $xy,x'y'\in F$ the subgraph induced by $x,y,x',y'$ has at least one edge not in $F$ different from $xy$ and $x'y'$.
  \end{itemize}
  \item[{\rm (iii)}] $S_F$ is a outer mutual-visibility set of $L(K_n)$ if and only if
    \begin{itemize}
    \item for any two not incident edges $uv,u'v'$ with $uv\in E(K_n)\setminus F$ and $u'v'\in F$ the subgraph induced by $u,v,u',v'$ has at least one edge not in $F$ different from $uv$ and $u'v'$, and
    \item for any two not incident edges $xy,x'y'\in F$ the subgraph induced by $x,y,x',y'$ has at least one edge not in $F$ different from $xy$ and $x'y'$.
  \end{itemize} 
\end{itemize}
\end{lemma}

By using Lemma~\ref{lem:TOD-L-K_n}, we can give the following conclusions on the (total, outer or dual) mutual-visibility number of $L(K_n)$.

\begin{proposition}
\label{prop:total-line-of-complete}
For any integer $n\ge 3$, $\mut(L(K_n))\ge n-1+\left\lfloor\frac{n-1}{2}\right\rfloor$.
\end{proposition}

\proof
Let $w\in V(K_n)$ and let $A=\{vw\in E(K_n)\,:\,v\in V(K_n)\setminus\{w\}\}$. Also, let $G'$ be the complete graph induced by $V(K_n)\setminus\{w\}$, and let $M$ be a maximum matching in $G'$. Now, consider the set of edges $F=A\cup M$ of $K_n$. Observe that $S_F$ satisfies the properties of Lemma \ref{lem:TOD-L-K_n} (i). Thus $S_F$ is a total mutual-visibility set of $L(K_n)$, and the bound follows since $|S_F| = n-1+\left\lfloor\frac{n-1}{2}\right\rfloor$.
\qed

By using a computer we have checked that the bound of Proposition~\ref{prop:total-line-of-complete} is tight for $n\in \{4,5,6,7\}$. On the other hand, the equality does not hold in general. For instance, the set $\{01,12,23,34,45,56,67,78,89,90, 04,19,26,38,57,79\}$ of vertices of $L(K_{10})$, or equivalently, edges of $K_{10}$, where we have taken $V(K_{10})=\{0,1,\dots,9\}$, is a total mutual-visibility set of $L(K_{10})$ of cardinality $16$. However, Proposition~\ref{prop:total-line-of-complete} only yields $\mut(L(K_{10}))\ge 13$.

The total mutual visibility of $L(K_n)$ has an interesting relation with the extension of the Tur\'an problem to forbidden generic graphs. 

\begin{definition}{\rm \cite[Page 479]{west-2021}}
The \emph{Tur\'an number} of a graph $H$, written $\ex(n;H)$, is the maximum number of edges in an $n$-vertex graph not containing $H$. 
\end{definition} 

\begin{theorem}\label{th:mut-Turan}
For any integer $n\ge 3$,  $\mut(L(K_n)) = \ex(n; C_4)$.
\end{theorem}
\begin{proof}
By Lemma~\ref{lem:TOD-L-K_n} (i), given a set $F\subseteq E(G)$, the set $S_F$ is a total mutual-visibility set of $L(K_n)$ if and only if for any two not incident edges
$uv, u'v' \in E(K_n)$, the subgraph induced by $u, v, u', v'$ has at least one edge not in $F$ different from $uv$ and $u'v'$. Consider any four vertices $u, v, u', v'$ of $K_n$. They induce a graph $G'=K_4$ with three pairs of not incident edges. Since for each pair at least one edge (not belonging to the pair) is not in $F$, it holds that at least two incident edges of $G'$ are not in $F$. Equivalently, this happens if and only if the edges of $F$ in $G'$ does not form a cycle $C_4$. Hence, in order to find a $\mut$-set in $L(K_n)$ we need to find a largest set of edges of $K_n$ that does not induce any $C_4$. By definition, its size is $\ex(n; C_4)$.\qed
\end{proof}

\begin{corollary}
 For any large enough integer $n$,  $\frac 1 2 (n^{3/2}-n^{4/3}) \le \mut(L(K_n)) \le \frac{1}{4} n(1+ \sqrt{4n-3})$.
\end{corollary}
\begin{proof}
    Given Theorem~\ref{th:mut-Turan} on the equivalence of $\mut(L(K_n))$ and $\ex(n; C_4)$, the upper bound was first proved by Reiman in~\cite{reiman-1958}. The lower bound (and a rediscovery of the upper bound) can be found in~\cite{brown-1966} 
    and~\cite{erdos-1966}. 
    \qed
\end{proof}

We next proceed with finding similar results as the above ones for the other two remaining mutual-visibility parameters (outer and dual).

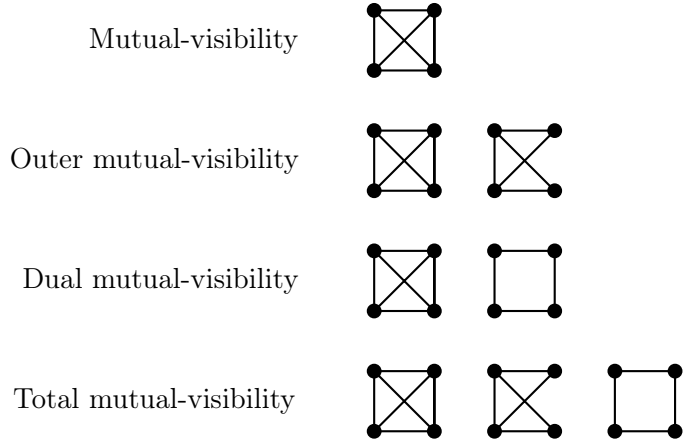
\begin{figure}[ht!]
\begin{center}
\begin{tikzpicture}[scale=0.4,style=thick,x=1cm,y=1cm]
\def\vr{6pt}

\begin{scope} 
\coordinate(a) at (0,14);
\coordinate(b) at (2,14);
\coordinate(c) at (2,12);
\coordinate(d) at (0,12);
\draw (a) -- (b) -- (c) -- (d) -- (a) -- (c) -- (b) -- (d);

\draw(a)[fill=black] circle(\vr);
\draw(b)[fill=black] circle(\vr);
\draw(c)[fill=black] circle(\vr);
\draw(d)[fill=black] circle(\vr);
\node at (-6,13) {Mutual-visibility};
\end{scope}
\begin{scope} 
\coordinate(a) at (0,10);
\coordinate(b) at (2,10);
\coordinate(c) at (2,8);
\coordinate(d) at (0,8);
\coordinate(e) at (4,10);
\coordinate(f) at (6,10);
\coordinate(g) at (6,8);
\coordinate(h) at (4,8);
\draw (a) -- (b) -- (c) -- (d) -- (a) -- (c) -- (b) -- (d);
\draw (e) -- (h) -- (g) -- (e) -- (f) -- (h);

\draw(a)[fill=black] circle(\vr);
\draw(b)[fill=black] circle(\vr);
\draw(c)[fill=black] circle(\vr);
\draw(d)[fill=black] circle(\vr);
\draw(e)[fill=black] circle(\vr);
\draw(f)[fill=black] circle(\vr);
\draw(g)[fill=black] circle(\vr);
\draw(h)[fill=black] circle(\vr);
\node at (-7.3,9) {Outer mutual-visibility};
\end{scope}
\begin{scope} 
\coordinate(a) at (0,6);
\coordinate(b) at (2,6);
\coordinate(c) at (2,4);
\coordinate(d) at (0,4);
\coordinate(e) at (4,6);
\coordinate(f) at (6,6);
\coordinate(g) at (6,4);
\coordinate(h) at (4,4);
\draw (a) -- (b) -- (c) -- (d) -- (a) -- (c) -- (b) -- (d);
\draw (e) -- (f) -- (g) -- (h) -- (e);

\draw(a)[fill=black] circle(\vr);
\draw(b)[fill=black] circle(\vr);
\draw(c)[fill=black] circle(\vr);
\draw(d)[fill=black] circle(\vr);
\draw(e)[fill=black] circle(\vr);
\draw(f)[fill=black] circle(\vr);
\draw(g)[fill=black] circle(\vr);
\draw(h)[fill=black] circle(\vr);
\node at (-7.1,5) {Dual mutual-visibility};
\end{scope}

\begin{scope} 
\coordinate(a) at (0,2);
\coordinate(b) at (2,2);
\coordinate(c) at (2,0);
\coordinate(d) at (0,0);
\coordinate(e) at (4,2);
\coordinate(f) at (6,2);
\coordinate(g) at (6,0);
\coordinate(h) at (4,0);
\coordinate(i) at (8,2);
\coordinate(l) at (10,2);
\coordinate(m) at (10,0);
\coordinate(n) at (8,0);
\draw (a) -- (b) -- (c) -- (d) -- (a) -- (c) -- (b) -- (d);
\draw (e) -- (h) -- (g) -- (e) -- (f) -- (h);
\draw (i) -- (l) -- (m) -- (n) -- (i);

\draw(a)[fill=black] circle(\vr);
\draw(b)[fill=black] circle(\vr);
\draw(c)[fill=black] circle(\vr);
\draw(d)[fill=black] circle(\vr);
\draw(e)[fill=black] circle(\vr);
\draw(f)[fill=black] circle(\vr);
\draw(g)[fill=black] circle(\vr);
\draw(h)[fill=black] circle(\vr);
\draw(i)[fill=black] circle(\vr);
\draw(l)[fill=black] circle(\vr);
\draw(m)[fill=black] circle(\vr);
\draw(n)[fill=black] circle(\vr);
\node at (-7.3,1) {Total mutual-visibility};
\end{scope}
\end{tikzpicture}
\caption{Forbidden induced subgraphs for $(K_n)_F$, where $F$ is such that $S_F$ is a (outer, dual, total) mutual visibility set of $L(K_n))$}
	\label{fig:forbidden}
\end{center}
\end{figure}

Theorems~\ref{th:mu-Turan} and~\ref{th:mut-Turan} provide us with a way to calculate the values of $\mu(L(K_n))$ and $\mut(L(K_n))$. They are based on the analysis of forbidden subgraphs for $(K_n)_F$ where $S_F$ is a mutual-visibility or a total mutual-visibility sets of $L(K_n))$, respectively.  By using Lemma~\ref{lem:TOD-L-K_n}, an analysis of the \emph{induced} forbidden subgraphs of $(K_n)_F$ for the (dual, outer, total) mutual-visibility sets $S_F$ of $L(K_n)$ shows that only three forbidden graphs are involved: $K_4$, $K_4^-$, and $C_4$ (see Figure~\ref{fig:forbidden}). As proved in Theorem~\ref{th:mu-Turan}, $K_4$ is the only forbidden subgraph of $(K_n)_F$ for any mutual-visibility sets $S_F$ of $L(K_n))$.
As for the total mutual-visibility, all the three induced graphs are forbidden, but, since $C_4$ is a subgraph of both $K_4$ and $K_4^-$, it is sufficient to forbid only this graph, and then $\mut(L(K_n))=\ex(n,C_4)$, as stated by Theorem~\ref{th:mut-Turan}. 

Similarly, for the outer mutual-visibility, the induced forbidden subgraphs are $K_4$ and $K_4^-$, and since $K_4^-$ is a subgraph of $K_4$, we have the following result.

\begin{theorem}\label{th:muo-Turan}
  For any integer $n\ge 3$, $\muo(L(K_n)) = \ex(n; K_4^-)$.
\end{theorem}

Based on this relationship above, and using the next known result, we are able to give the exact value of $\muo(L(K_n))$.

\begin{theorem}{\em \cite{Simonovits}}
If $F$ has chromatic number $k$ and a critical edge, and $n$ is large enough, then $\ex(n, F) = |E(T(n,k-1))|$. Moreover, $T(n,k-1)$ is the unique extremal graph.
\end{theorem}

Since the graph $K_4^-$ has chromatic number $3$ and a critical edge, we deduce that the edges of $K_n$ that form an outer mutual-visibility set of the largest cardinality in $L(K_n)$, together with the vertices in such edges, form a graph isomorphic to the Tur\'an graph $T(n,2)$. Recall that $T(n,2)$ is the bipartite graph of order $n$ with partite sets of cardinality $\lceil n/2\rceil$ or $\lfloor n/2\rfloor$. Thus, the following result holds.

\begin{corollary}
  For any large enough integer $n$, $\muo(L(K_n)) = \lceil \frac{n}{2}\rceil\cdot\lfloor \frac{n}{2}\rfloor$.
\end{corollary}

Now, for the dual mutual-visibility, the two induced forbidden subgraphs are $K_4$ and 
$C_4$. Then the following result holds.

\begin{theorem}\label{th:mud-Turan}
  Let  $F \subseteq E(K_n)$. Then $S_F \subseteq V (L(K_n))$ is a dual mutual-visibility set of
$L(K_n)$  if and only if $(K_n)_F$ is a $(K_4,C_4)$-free graph.
\end{theorem}

In contrast with the cases which appear along with standard, total and outer mutual-visibility sets, there is a lack (to the best of our knowledge) of results concerning the largest number of edges in a  $(K_4,C_4)$-free graph of order $n$. This made that the result above for the dual mutual-visibility sets does not lead to a bound or formulae for the dual mutual-visibility number of $L(K_n)$. By computer checking, we only know that for $n\in[10]$ the largest  graphs have $0,1,3,5,7,10,12,15,18,21$ edges, respectively. On the positive side, this means that it is worthy of considering studying this problem independently.

\medskip
We next continue with the line graphs of a complete bipartite graph $K_{m,n}$, $m, n\ge 2$. Then we recall that Palmer~\cite{palmer-1973} proved that the line graph of a connected graph $G$ is a nontrivial Cartesian product if and only if $G = K_{n,m}$, $n,m\ge 2$, see~\cite[Proposition 1.2]{ikr-2008}. So, $L(K_{m,n}) \cong K_m\cp K_n$. It is already known from~\cite{Cicerone-2023} that $\mu(K_m\cp K_n) = z(m,n;2,2)$, where $z(m, n; 2, 2)$ is the Zarankiewicz number, that can also be seen as the maximum number of edges in a complete bipartite graph $K_{m,n}$ that has no $4$-cycle. We now state that the same conclusion can be also obtained by using Lemma~\ref{lem:MVS-line}. The proof of it runs along the lines of the proof of Theorem \ref{th:mu-Turan}, and thus it is left to the reader.

\begin{theorem}
\label{thm:LKmn}
Let $n,m\ge 2$ and $F\subseteq E(K_{m,n})$. Then $S_F\subseteq V(L(K_{m,n}))$ is a $\mu$-set of $L(K_{m,n})$ if and only if $S_F$ is a largest set of vertices of $L(K_{m,n})$ such that $(K_{m,n})_F$ contains no $4$-cycle.
\end{theorem}

The following consequence of Theorem~\ref{thm:LKmn} then follows from the above-mentioned observations from~\cite{Cicerone-2023}.

\begin{corollary}
\label{cor:line-of-complete-bipartite}
For any two integers $n,m\ge 2$, $\mu(L(K_{m,n}))=z(m, n; 2, 2)$.
\end{corollary}

We close this section by again using the fact that the line graph of a complete bipartite graph $K_{n,m}$ is isomorphic to $K_n\cp K_m$. Hence, a result from \cite{tian-2023+}, and Theorem \ref{thm:hamming-outer-dual} lead to the following consequence.

\begin{corollary}
\label{cor:line-of-complete-bipartite-others}
For any two integers $n,m\ge 2$, $\mut(L(K_{m,n}))=\max\{n,m\}$, $\mud(L(K_{m,n}))=n+m-1$, and $\muo(L(K_{m,n}))=n+m-2$.
\end{corollary}

\section{Visibility in cographs}
\label{sec:cographs}


In this section, we consider the mutual-visibility in cographs. In the main result we prove that if $G$ is a cograph, then either $\mu(G) = \mut(G)$ ($= \muo(G) = \mud(G)$), or $\mu(G) = \mud(G) = n(G)-1$ and $\mut(G) = \muo(G) =n(G)-2$.

A cograph is a graph all of whose connected induced subgraphs have diameter at most 2. Moreover, each cograph can also be built up from a single vertex by adding a sequence of \emph{twins}. Two vertices $u$ and $v$ of $G$ are twins if $N_G(u)=N_G(v)$; in particular, they are \emph{false twins} if they are twins and $uv\notin E(G)$, and they are \emph{true twins} if they are twins and $uv\in E(G)$. 
We use some additional notation here. Given a graph $G$ and $v\in V(G)$, $G-v$ denotes the subgraph of $G$ induced by $V(G)\setminus \{v\}$. A graph $G$ such that  $\mu(G)=\mut(G)$ is called a {\em $(\mu, \mut)$-graph}. 

We start by recalling the following characterizations from~\cite{DiStefano-2022} and~\cite{Cicerone-2022+}. 

\begin{lemma}\label{lem:old_n-1}
{\rm \cite[Lemma 4.8]{DiStefano-2022} }
Given a graph $G$, then $\mu(G)\ge n(G)-1$ if and only if there exists a vertex $v$ in $G$ adjacent to each vertex $u$ in $G-v$ such that $\deg_{G-v} (u) < n(G) - 2$.
\end{lemma}

In~\cite{Cicerone-2022+}, any vertex $v$ of $G$ fulfilling the condition in the above lemma was called \emph{enabling}.

\begin{proposition}\label{prop:old_cographs}
{\rm \cite[Proposition 3.5]{Cicerone-2022+} }
A cograph $G$ is a $(\mu, \mut)$-graph if and only if it has a universal vertex or no enabling vertices.
\end{proposition}

The following definition aims to reformulate the previous characterizations in terms of graph structure.

\begin{definition}
A \emph{big-$\mu$} graph is any graph $G$ defined as $G = (K_1 \cup K_t) + H$, where $K_1$, $K_t$, and $H$ are three distinct graphs such that $t\ge 0$ (i.e., $K_t$ can be an empty graph). 
\end{definition}

From this definition, it follows that each non-trivial clique is a big-$\mu$ graph (it is sufficient to take $K_t$ empty and $H$ as a clique). Consequently, observe that if $G$ is a big-$\mu$ graph, then $\mu(G)=n(G)$ when $K_t$ is empty and $H$ isomorphic to a clique, and $\mu(G)=n(G)-1$ otherwise. This observation explains the term \emph{big-$\mu$}. 

\smallskip
The two characterizations recalled above will be reformulated by using the following lemma.

\begin{lemma}\label{lem:big_equiv_enabling}
Let $G$ be an arbitrary graph. Then $G$ is a \emph{big-$\mu$ graph} if and only if there exists a vertex $v$ in $G$ adjacent to each vertex $u$ in $G-v$ such that $\deg_{G-v} (u) < n(G) - 2$.
\end{lemma}

\begin{proof}
$(\Rightarrow)$
Assume that $G$ is a big-$\mu$ graph, that is there exists three distinct graphs $K_1$, $K_t$, and $H$  such that $G= (K_1\cup K_t) + H$. Let $v$ be the vertex forming the graph $K_1$.  Since the vertices $u$ of $H$ are the only vertices such that $\deg_{G-v} (u) < n(G) - 2$, and since $v$ is adjacent to all of them, the thesis follows.

$(\Leftarrow)$
Assume now there exists a vertex $v$ of $G$ being adjacent to each vertex $u$ in $G-v$ such that $\deg_{G-v} (u) < n(G) - 2$. To show that $G$ is a big-$\mu$ graph, take $K_1$ formed by $v$, $G'$ as the graph induced by $N(v)$, and $G''$ as the graph induced by $V(G)\setminus N[v]$.  
Let $u''$ be a vertex of $G''$. Since $u''$ is not adjacent to $v$ in $G$, it holds that $\deg_{G-v} (u'') \ge n(G) - 1$. This implies that $G''$ is a clique, and there exists an edge in $G$ between each pair of vertices $u'\in V(G')$ and $u''\in V(G'')$.
Hence $G=(K_1\cup G') + G''$, where $G'$ is an arbitrary graph and $G''$ is a clique. 
\qed 
\end{proof}

\begin{corollary}\label{cor:new_n-1}
Let $G$ be an arbitrary graph. Then $\mu(G)\ge n(G)-1$ if and only if $G$ is a big-$\mu$ graph.
\end{corollary}

\begin{proof}
This is an immediate consequence of Lemma~\ref{lem:old_n-1} and Lemma~\ref{lem:big_equiv_enabling}.
\qed 
\end{proof}

\begin{corollary}\label{cor:new_cographs}
Let $G$ be a cograph. Then $\mu(G) > \mut(G)$ if and only if $G$ is a big-$\mu$ graph $G=(K_1\cup K_t)+H$ 
with no universal vertices.
\end{corollary}
\begin{proof}
$(\Rightarrow)$
From Proposition~\ref{prop:old_cographs} we have $\mu(G) > \mut(G)$ if and only if $G$ has an enabling vertex $v$ and $G$ has no universal vertices. By Lemma~\ref{lem:big_equiv_enabling}, we have that $G=(K_1\cup K_t)+H$ is a big-$\mu$ graph.

$(\Leftarrow)$
As $G=(K_1 \cup K_t) + H$ has no universal vertices, $G$ is not a clique and then $\mu(G)< n(G)$. Hence, by Corollary~\ref{cor:new_n-1}, $\mu(G)= n(G)-1$.
We claim that $\mut(G)< n(G)-1$.
Assume, on the contrary, that $\mut(G)= n(G)-1$.
Let $S$ be a $\mut$-set of $G$ and let $u$ be the only vertex of $V(G)$ not in $S$ and let $v$ be the only vertex in $K_1$.
Since $G$ has no universal vertices we deduce, (1) $K_t$ is not empty and (2) $H$ has at least two not adjacent vertices $x$ and $y$. Then, by (1),  $u$ cannot be $v$, otherwise $u$ is not in mutual visibility with any vertex in $K_t$. Moreover, $u$ cannot be a vertex of $K_t$, otherwise $u$ is not in mutual visibility with $v$. By (2), $u$ cannot be a vertex of $H$, otherwise $x$ and $y$ are not in mutual visibility.
Hence $\mut(G)< n(G)-1$.
%
\qed
\end{proof}

This corollary implies that the smallest cograph $G$ which is not a $(\mu,\mut)$-graph corresponds to the cycle $C_4 = (K_1\cup K_t) + H$, with $t=1$ and $H=K_1\cup K_1$. If $v$ is the vertex forming $K_1$, $h_1,h_2$ are the vertices forming $H$, and $k$ is the unique vertex of $K_t$, it can be observed that each cograph which is not a $(\mu,\mut)$-graph can be obtained from this initial cycle by applying no split operations to $v$, any possible split operation to $h_1$ and $h_2$ and  only true-twin operations to $k$.

\smallskip
In~\cite{DiStefano-2022} it is shown that $\mu(G)\ge n(G)-2$ for each cograph $G$, and that the exact value of $\mu(G)$ can be computed in polynomial time. The following statement extends the analysis to the other visibility parameters.

\begin{theorem}\label{thm:values_cographs}
If $G$ is a cograph, then either $\mu(G) = \mut(G)$, or $\mu(G) = \mud(G) = n(G)-1$ and $\mut(G) = \muo(G) =n(G)-2$.
\end{theorem}
\begin{proof}
If $\mu(G) \not = \mut(G)$ then, by Corollary~\ref{cor:new_cographs}, $G$ is a big-$\mu$ cograph $G=(K_1\cup K_t)+H$  without universal vertices. By the proof of the same corollary, we have $\mu(G) = n(G) - 1$. Concerning $\mud(G) = n(G)-1$, it easily follows by observing that $V(G)\setminus \{u\}$, where $u$ is the unique vertex of $K_1$, is a dual mutual-visibility set of $G$. 

By using again the proof of  Corollary~\ref{cor:new_cographs}, we have $\mut(G) < n(G) - 1$. Notice that the same arguments can be used to show that $\muo(G) < n(G) - 1$ also holds. 
Hence, to conclude the proof, it is sufficient to show there exists a total mutual-visibility set of $G$ with $n(G)-2$ elements. 
To this end, let $S=V(G)\setminus \{u,v\}$ where $u$ is the unique vertex of $K_1$ and $v$ is a vertex in $H$. Let us show that $S$ is a total mutual-visibility set of $G$. 
The vertices in $S$ are in mutual visibility, since any two of them are adjacent or they are adjacent vertices of $u$.
The vertices $u$ and $v$ are adjacent. Vertex $u$ is in mutual visibility with all the vertices in $S$, as $u$ is adjacent to them or in mutual visibility through vertex $v$. Analogously, vertex $v$ is in mutual visibility with all the vertices in $S$ as $v$ is adjacent to them or in mutual visibility through vertex $u$. Then, $S$ is a 
total mutual-visibility set of $G$ with $n(G)-2$ vertices.
\qed
\end{proof}

A well-known superclass of cographs is that formed by \emph{distance-hereditary graphs}. In fact, these graphs can be generated by using true twins, false twins, and pendant vertices. Concerning the problem of characterizing all the distance-hereditary graphs graphs $G$ in which  $\mu(G) > \mut(G)$, we conjecture that the following holds:

\begin{conjecture}
If $G$ is a distance-hereditary graph but not a big-$\mu$ cograph without universal vertices, then $\mu(G) = \mut(G)$. 
\end{conjecture}

We must remark that this conjecture is also supported by numerous computer-assisted simulations.

\section{Non-trivial diameter-two graphs of minimum size}
\label{sec:of-minimum-size}

Let $G$ be a diameter-two graph with no universal vertex. Then Erd\H{o}s and R\'{e}nyi proved that $m(G)\ge 2n(G) - 5$. More than two decades later, Henning and Southey characterized the graphs which achieve the bound. In this section we determine $\mu, \muo, \mud$, and $\mut$ for these extremal diameter-two graphs. 

Let us restate the mentioned classical result of Erd\H{o}s and R\'{e}nyi on the minimum size of a diameter-two graph with no universal vertex.

\begin{theorem}\label{thm:erdos-renyi}
{\rm \cite{erdos-1962} } 
If $G$ is a diameter-two graph with no universal vertex, then $m(G)\ge 2n(G) - 5$.
\end{theorem}

Let $P$ be the Petersen graph and note that $P$ attaines the bound of Theorem~\ref{thm:erdos-renyi}. It is already known that $\mu(P) = 6$, see~\cite{Cicerone-2023}, and that $\mut(P) = 0$, see~\cite{tian-2023+}. By a case analysis, we also get that $\mud(P) = 0$. On the other hand, Lemma~\ref{lem:girth5} implies that $\muo(P)\le 4$, and it can be easily checked that an independent set of $P$ of cardinality $4$ is an outer mutual-visibility set. In summary, 
$$\mut(P) = \mud(P) = 0, \quad \muo(P) = 4,\quad \mu(P) = 6\,.$$

\begin{figure}[ht!]
\begin{center}
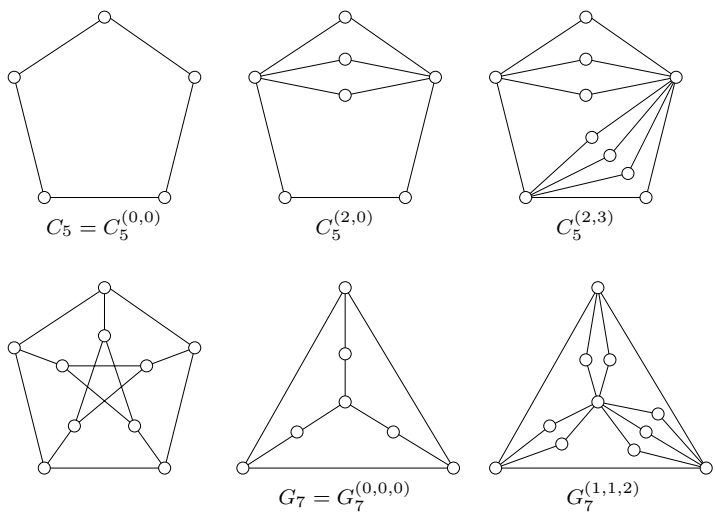
\end{center}
\caption{Some examples of graphs from the family $\mathcal{G}$.}
\label{fig:family}
\end{figure}

Let $G_7$ (cf. Fig.~\ref{fig:family}) be the graph obtained from the cycle $C_3$ by adding a pendant edge to each vertex of the cycle and then adding a new vertex and joining it to the three degree-one vertices. In~\cite{henning-2015}, the following family $\mathcal{G}$ of graphs has been defined: 
\begin{enumerate}
\item[(i)] $\mathcal{G}$ contains $C_5$, $G_7$, and the Petersen graph; and 
\item[(ii)] $\mathcal{G}$ is closed under degree-2 vertex duplication (cf.\ Fig.~\ref{fig:family}).
\end{enumerate}

The graphs that achieve equality in the bound of Theorem~\ref{thm:erdos-renyi} are characterized as follows. 

\begin{theorem}
{\rm \cite{henning-2015} } 
If $G$ is a diameter-two graph of order $n$ and size $m$ with no universal vertex, then $m = 2n - 5$ if and only if $G\in \mathcal{G}$.
\end{theorem}

In what follows, we compute $\tau(G)$ for each $G\in \mathcal{G}$ and for each $\tau\in \{\mu,\muo,\mud,\mut\}$. To this aim, denote by $C_5^{(i,j)}$ the graph obtained from the cycle $C_5$ with one vertex duplicated $i\ge 0$ times, and another vertex duplicated $j\ge 0$ times. For instance, in Fig.~\ref{fig:family} the graphs $C_5 = C_5^{(0,0)}$, $C_5^{(2,0)}$, and $C_5^{(2,3)}$ are shown.

\begin{figure}[ht!]
\begin{center}
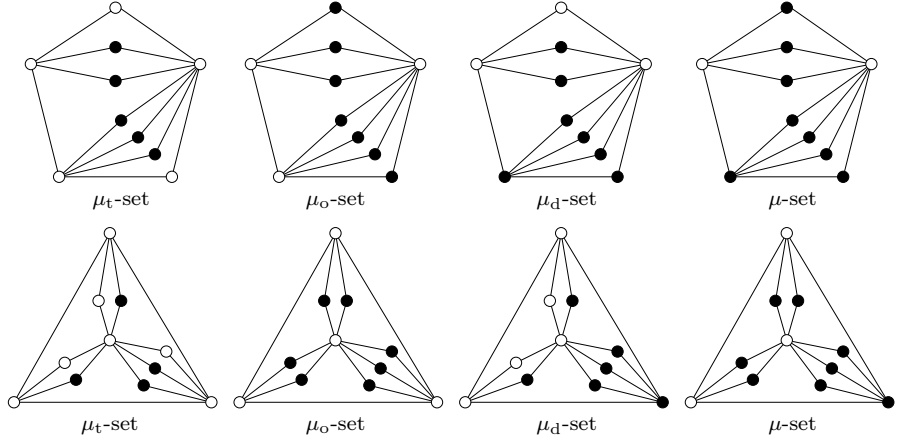
\end{center}
\caption{Variations of mutual-visibility sets in some graphs of the family $\mathcal{G}$.
}
\label{fig:family-sol}
\end{figure}

\begin{lemma}\label{lem:mu_subgraph}
If $X$ is a  mutual-visibility  (resp. outer, dual, or total mutual visibility) set of a graph $G$ and $x\in X$, then $X\setminus \{x\}$ is a  mutual-visibility  (resp. outer, dual, or total mutual visibility) set of $G - x$.
\end{lemma}

\begin{proof}
Let $X$ be a mutual-visibility set. For any two vertices $u,v$ in $X\setminus \{x\}$, there is at least one $u,v$-shortest path $P$ not passing through $x$. Then the removal of $x$ from $G$ does not destroy $P$ in  $G - x$. Hence, for the generality of $u,v$, the set $X\setminus \{x\}$ is a mutual-visibility set of $G$. Similar arguments work for the outer, dual or total mutual-visibility cases.
\qed
\end{proof}

\begin{proposition}\label{prop:C5}
If $i, j\ge 0$, then 
\begin{itemize}
\item 
$\mut(C_5^{(i,j)}) = i+j$, 
\item 
$\muo(C_5^{(i,j)}) = \mud(C_5^{(i,j)}) = i+j+2$, 
\item $\mu(C_5^{(i,j)}) = i+j+3$.
\end{itemize}
\end{proposition}
\begin{proof}
In~\cite{variety-2023} it is shown that the four formulae are correct when $i=j=0$. (Notice, however, that although it is $\muo(C_5^{(0,0}) = \mud(C_5^{(0,0)}) = 2$, two vertices in a dual mutual-visibility are adjacent, whereas two vertices in an outer mutual-visibility are not.) Fig.~\ref{fig:family-sol} shows a visibility set for each of the four variants. In each case, it is clear that each time a degree-2 vertex duplication is made, the new vertex can be included in the visibility-set and hence the corresponding parameter is increased by one. 

Concerning the optimality, consider now all the four cases of mutual-visibility sets for  $G=C_5^{(i,j)}$. We prove by induction that $\mut(G) = i+j$, $\muo(G) = i+j+2$, $\mut(G) = i+j+2$, $\mut(G) = i+j+3$. As already observed, the statement holds for the initial case in which $i=j=0$. Assume it holds for $k=i+j>0$ and consider the case $G= C_5^{(i,j)}$ obtained with $k+1$ degree-2 vertex duplications. Assume, by contradiction, that there exists a $\mut$-set ($\muo$-set, $\mud$-set, $\mu$-set) $X$ of $G$ with $|X| > n(G)-5$ ($|X| > n(G)-3$, $|X| > n(G)-3$, $|X| > n(G)-2$). Then we can consider any vertex $v\in X$ and set $X' = X\setminus \{v\}$ and $G' = G-v$. By Lemma~\ref{lem:mu_subgraph}, $X'$ is a total mutual-visibility (outer mutual-visibility, dual mutual-visibility, mutual-visibility) set of $G'$, with a size larger than that assumed by induction. 
\qed
\end{proof}
 
%


It is worth to remark that the case concerning $\mu(G)$ in Proposition~\ref{prop:C5} can be also proved by simply observing that the mutual-visibility set of $G$, shown in  Fig.~\ref{fig:family-sol}, contains $i+j+3 = n(G)-2$ elements, and that Corollary~\ref{cor:new_n-1} implies that this set is indeed a $\mu$-set, since $G$ is not a big-$\mu$ graph. 


Similarly as above, denote by $G_7^{(i,j,k)}$ the graph obtained from $G_7$ when the three degree-2 vertices have been respectively duplicated $i$, $j$, and $k$ times. For instance, Fig.~\ref{fig:family} shows the graphs $G_7 = G_7^{(0,0,0)}$ and $G_7^{(1,1,2)}$.

\begin{proposition}\label{prop:G7}
For the graph $G_7^{(i,j,k)}$ we have the following formulae: 
\begin{itemize}
\item 
$\mut(G_7^{(i,j,k)}) = i+j+k$, 
\item 
$\muo(G_7^{(i,j,k)}) = i+j+k+3$, 
\item 
$\mud(G_7^{(i,j,k)}) = \begin{cases}
  3;  & i+j+k = 0, \\
  i+j+k+2; &  i+j+k\ge 1.
\end{cases}$
\item 
$\mu(G_7^{(i,j,k)}) = i+j+k+4$.
\end{itemize}
\end{proposition}
\begin{proof}
The statement can be proved by using the same inductive approach used in the proof of Proposition~\ref{prop:C5}. The correctness of the base cases can be easily verified thanks to the limited size of the graph. We just remark that for the dual mutual-visibility case, when  $i=j=k=0$, the $\mud$-set is composed of the three vertices of the $C_3$ cycle; when the first degree-2 vertex duplication is made (i.e., $i=1$ and $j=k=0$), a $\mud$-set can be identified by selecting the two ``twin'' degree-2 vertices along with their adjacent vertex in the $C_3$ cycle. Extending this graph further leads to identifying the $\mud$-set as represented in Fig.~\ref{fig:family-sol}.
\qed
\end{proof}

\section{Concluding remarks}
\label{sec:conclusion}

In this paper we considered graphs of diameter two and their values for the (classic, total, dual and outer) mutual-visibility parameters. We next comment some possible open questions that might be of interest to continue exploring.

\begin{itemize}
    \item The class of graphs of diameter two is very wide. We have studied here a few of them, but some other non-trivial classes might be of interest as well. Among them, we remark the Kneser graph $K(n, 2)$, which is indeed the complement of $L(K_n)$, and the line graphs of complete multipartite graphs (with at least three partite sets). For this latter ones, since we do not have an exact solution for $\mu(L(K_{m,n}))$ (and it seems to be beyond the reach of existing methods), we cannot, of course, expect an exact result for the mutual-visibility number of general complete multipartite graphs. Suppose $F\subseteq E(K_{n_1,\ldots, n_k})$, $k\ge 3$, $n_1, \ldots, n_k\ge 2$, is a $\mu$-set of $L(K_{n_1,\ldots, n_k})$. As $\diam(L(K_{n_1,\ldots, n_k})) = 2$, we can again apply Lemma~\ref{lem:MVS-line} to see that if $uv,u'v'\in F$ are independent edges, then there is an edge $xy$ incident with both $uv$ and $u'v'$ such that $xy\notin F$. If the edges $uv$ and $u'v'$ are only between two of the multipartite sets, then in $(K_{n_1,\ldots, n_k})_F$ we have $C_4$ as a forbidden subgraph. If the edges $uv$ and $u'v'$ lie in three multipartite sets, then in $(K_{n_1,\ldots, n_k})_F$ we have $K_4^-$ as a forbidden subgraph, while if the end vertices of $uv$ and $u'v'$ lie in four multipartite sets, then in $(K_{n_1,\ldots, n_k})_F$ we have $K_4$ as a forbidden subgraph. However, all these facts are not exactly related to each other, since the situations are somehow not comparable. Consequently, it would be interest to continue the study of the (dual, outer, total) mutual-visibility number of these line graphs. 
    \item Theorem \ref{thm:hamming-outer-dual} completes the studies on the mutual-visibility variants of $2$-dimensional Hamming graphs (those of diameter two). For higher dimensions, the total version was studied in \cite{Bujtas}, and the problem seems to be very challenging due to its connection with some Tur\'an type problems in hypergraphs. This makes natural to consider the remaining variants (dual and outer) for Hamming graphs of higher dimension.  
    \item Based on the fact that finding the value of any of the studied mutual-visibility parameters of graphs of diameter two seems to be a hard task, we consider the following question of interest. Which is the computational complexity of computing the (outer, dual, total and classical) mutual visibility number of graphs of diameter two?
    \item In connection with Theorem \ref{th:mud-Turan}, as we already mentioned, it seems there is a lack of results concerning the largest number of edges in a  $(K_4,C_4)$-free graph of order $n$. Based on this fact, it might be of interest to separately study this problem from a combinatorial point of view. A consequence of such study will clearly give some knowledge on the dual mutual-visibility number of $L(K_n)$.
\end{itemize}

\section*{Acknowledgments}
S.\ Cicerone and G.\ Di Stefano were partially supported by the European project  ``Geospatial based Environment for Optimisation Systems Addressing Fire Emergencies'' (GEO-SAFE), contract no. H2020-691161, and by the Italian National Group for Scientific Computation (GNCS-INdAM).
S.\ Klav\v{z}ar was partially supported by the Slovenian Research Agency (ARRS) under the grants P1-0297, J1-2452, and N1-0285.
I.\ G.\ Yero has been partially supported by the Spanish Ministry of Science and Innovation through the grant PID2019-105824GB-I00. Moreover, this investigation was
completed while this author (I.G.\ Yero) was making a temporary stay at ``Universitat Rovira i Virgili'' supported by the program ``Ayudas para la recualificaci\'on del sistema universitario espa\~{n}ol para 2021-2023, en el marco del Real Decreto 289/2021, de 20 de abril de 2021''.

\section*{Author contributions statement} All authors contributed equally to this work.

\section*{Conflicts of interest} The authors declare no conflict of interest.

\section*{Data availability} No data was used in this investigation.


\end{document}